\documentclass{amsart}
\usepackage{bm}
\usepackage{amsthm}
\usepackage[usenames,dvipsnames]{pstricks}
\usepackage{epsfig}
\usepackage{pst-grad} 
\usepackage{pst-plot} 
\usepackage[small,nohug]{diagrams}

\newtheorem{theorem}{Theorem}[section]

\newtheorem{lemma}[theorem]{Lemma}
\newtheorem{proposition}[theorem]{Proposition}

\theoremstyle{remark}
\newtheorem{remark}[theorem]{Remark}

\theoremstyle{definition}
\newtheorem{definition}[theorem]{Definition}
\newtheorem{example}[theorem]{Example}

\newcommand{\ilim}{\varprojlim}

\newcommand{\R}{\mathrm{R}}

\newcommand{\N}{\mathrm{N}}
\newcommand{\cov}{\mathrm{cov}}

\newcommand{\St}{\mathrm{St}}

\begin{document}

\title{A new version of homotopical Hausdorff}
\author{B.~LaBuz}
\address{Saint Francis University, Loretto, PA 15940}
\email{blabuz@@francis.edu}
\date{}

\begin{abstract}
It is known that shape injectivity implies homotopical Hausdorff and that the converse does not hold, even if the space is required to be a Peano continuum. This paper gives an alternative definition of homotopical Hausdorff inspired by a new topology on the set of fixed endpoint homotopy classes of paths. This version is equivalent to shape injectivity for Peano spaces. 
\end{abstract}

\maketitle
\tableofcontents

\section{Introduction}

Homotopical Hausdorff is a homotopical criterion that detects if a Hausdorff space $X$ has $\widetilde X$ Hausdorff where $\widetilde X$ is the set of fixed endpoint homotopy classes of paths in $X$ starting at some basepoint with a standard topology (see Definition \ref{BasicTopology}). Homoropical Hausdorff was discussed in \cite{Zas} and \cite{BogSie} and given its present name in \cite{CanCon}. In \cite{FisZas} the authors show that shape injectivity implies homotopical Hausdorff. A space $X$ is shape injective if the homomorphism $\pi_1(X)\to\check\pi_1(X)$ is injective. Spaces that are known to be shape injective include one dimensional Hausdorff compacta \cite{EdaKaw} and subsets of closed surfaces \cite{FisZas2}. Thus these spaces have nice spaces of path homotopy classes.

The authors in \cite{Conetal} give examples of two Peano continua that are homotopically Hausdorff but not shape injective. The present paper notes that if the definition of homotopical Hausdorff is modified to reflect a new topology on $\widetilde X$ (Definition \ref{SmallLoop}) then the two concepts are equivalent for Peano spaces (Theorem \ref{MainTheorem}).

Section \ref{ShapeInject} gives a treatment of shape injectivity that mirrors the theory of generalized paths in \cite{Rips}. This viewpoint relates paths in a space $X$ to chains of points in $X$. It is quite geometric and lends itself to the proof of Theorem \ref{MainTheorem}.

\section{Shape injectivity}\label{ShapeInject}

Generalized paths and the uniform shape group are defined for uniform spaces in \cite{Rips}. We introduce an analogous construction for all topological spaces. Let us first recall the definition of the classical shape group.

We will consider only normal open covers.   Recall an open cover $\mathcal U$ of $X$ is normal if it admits a partition of unity $\{\phi_U:X\to [0,1]\}_{U\in\mathcal U}$ with $\phi_U^{-1}(0,1]\subset U$ for each $U\in\mathcal U$. We say that the partition of unity is subordinate to $\mathcal U$. The partition of unity can be chosen to be locally finite. The nerve $\N(\mathcal U)$ of $\mathcal U$ is the simplicial complex whose vertices are elements of $\mathcal U$ and whose simplices are finite subsets $\mathcal A\subset \mathcal U$ such that the intersection of the elements of $\mathcal A$ is nonempty. 

Fix an $x_0\in X$ and for each normal cover $\mathcal U$ fix a $U_0\in\mathcal U$ with $x_0\in U_0$. Given two normal covers $\mathcal U$ and $\mathcal V$, define $\mathcal U\leq \mathcal V$ if $(\mathcal V,V_0)$ refines $(\mathcal U,U_0)$, that is, $\mathcal V$ refines $\mathcal U$ and $V_0\subset U_0$. In this case choose a bonding map $p:\N(\mathcal V)\to \N(\mathcal U)$ such that each $V\in \mathcal V$ gets mapped to a $U\in\mathcal U$ with $V\subset U$, making sure to send $V_0$ to $U_0$.  The shape group $\check\pi_1(X,x_0)$ is the inverse limit $\ilim \pi_1(\N(\mathcal U),U_0)$. \cite{MS}

Given a cover $\mathcal U$ of $X$ define a $\mathcal U$-chain to be a finite list $x_1,\ldots x_n$ of points in $X$ such that for each $i<n$, $x_i,x_{i+1}\in U$ for some $U\in \mathcal U$. Let $\R(X,\mathcal U)$ be the simplicial complex whose vertices are points in $X$ and $A\subset X$ is a simplex if it is a finite subset of some $U\in\mathcal U$. It is the Rips complex of $X$ with respect to $\mathcal U$.

We identify a $\mathcal U$-chain $x_1,\ldots,x_n$ with the concatenation of the edge paths $[x_1,x_2],\ldots,[x_{n-1},x_n]$ in $\R(X,\mathcal U)$. We define two $\mathcal U$-chains to be $\mathcal U$-homotopic if the corresponding paths in $\R(X,\mathcal U)$ are fixed endpoint homotopic. This homotopy can be chosen to be simplicial. Thus two $\mathcal U$-chains are $\mathcal U$-homotopic if one can move from one to the other by a finite sequence of vertex additions and deletions.

Define a generalized path in $X$ to be a collection $\alpha=\{[\alpha_{\mathcal U}]_{\mathcal U}\}$ of equivalence classes of $\mathcal U$-chains in $X$, where $\mathcal U$ runs over all normal open covers of $X$, such that if $\mathcal V$ refines $\mathcal U$, $\alpha_{\mathcal V}$ is $\mathcal U$-homotopic to $\alpha_{\mathcal U}$. We define the covering shape group $\check\pi_1^{\cov}(X,x_0)$ to be the group of generalized loops in $X$ based at $x_0$ under the operation of concatenation. It is isomorphic to $\ilim \pi_1(R(X,\mathcal U),x_0)$.

We will show that $\check\pi_1^{\cov}(X,x_0)$ is isomorphic to the classical shape group. In order to do so, let us recall the following definition. Given an open cover $\mathcal U$ of $X$, the star of a point $x\in X$ in $\mathcal U$ is the union of all $U\in \mathcal U$ containing $x$. We say that a cover $\mathcal V$ is a star refinement of a cover $\mathcal U$ if the cover $\{\St(x,\mathcal V):x\in X\}$ refines $\mathcal U$. Any normal open cover has a normal star refinement. \cite[Proposition 5.3]{PU}.

It is more convenient to use the following notion of a star of a cover. Given an open cover $\mathcal U$ and a $U\in \mathcal U$, let $\St U$ be the union of all $V\in\mathcal U$ that meet $U$. Let $\St\mathcal U$ be the set of all $\St U$ for $U\in\mathcal U$. Notice the similarity between the open set $\St U$ in $X$ and the open star $\St U$ of the vertex $U$ in $\N(\mathcal U)$. They are both defined in terms of all $V\in\mathcal U$ that meet $U$.

\begin{lemma}
Suppose $\mathcal W$ is a star refinement of $\mathcal V$ and that $\mathcal V$ is a star refinement of $\mathcal U$. Then $\St\mathcal W$ refines $\mathcal U$.
\end{lemma}

\proof
Given $W\in\mathcal W$, let $x\in W$. We will show that $\St(W,\mathcal W)\subset \St(x,\mathcal V)$. Suppose $y\in\St(W,\mathcal W)$. Then $y\in W'$ where $W'\in\mathcal W$ meets $W$, say at a point $z$. Then $x,y\in\St(z,\mathcal W)$ which is contained in some $V\in\mathcal V$. Thus $y\in\St(x,\mathcal V)$.
\endproof

\begin{proposition}
$\check\pi_1(X,x_0)$ is isomorphic to $\check\pi_1^{\cov}(X,x_0)$.
\end{proposition}

\proof
Fix a basepoint $x_0\in X$ and for each normal cover $\mathcal U$ of $X$, fix a ``basepoint'' $U_0\in\mathcal U$ with $x_0\in U_0$. Define a pointed map $(\R(X,\mathcal U),x_0)\to(\N(\St\mathcal U),\St U_0)$ to send a vertex $x\in X$ to $\St U$ for some $U\in\mathcal U$ with $x\in U$. Note we can assume $\St U_0$ is the basepoint of $\St \mathcal U$ since any other $\St U$ that contains $x_0$ meets $\St U_0$. Let us see that this map on vertices extends to a simplicial map. Suppose $[x_1,\ldots,x_n]\in \R(X,\mathcal U)$. If $x_i\mapsto \St U_i$, then $x_1\in\St U_i$ for each $i\leq n$ so $[\St U_1,\ldots,\St U_n]\in\N(\St\mathcal U)$.

Now define a pointed map $(\N(\mathcal U),U_0)\to (\R(X,\St\mathcal U),x_0)$ to send a vertex $U\in\mathcal U$ to some $x$ with $x\in U$. Let us see that this map extends to a simplicial map. Suppose $[U_1,\ldots,U_n]\in\N(\mathcal U)$. If $U_i\mapsto x_i$, then $x_1,\ldots,x_n\in\St U_1$ so $[x_1,\ldots,x_n]\in\R(X,\St\mathcal U)$.

These maps induce homomorphisms $\pi_1(\R(X,\mathcal U),x_0)\to \pi_1(\N(\St\mathcal U),\St U_0)$ and $\pi_1(\N(\mathcal U),U_0)\to \pi_1(\R(X,\St\mathcal U),x_0)$. By the lemma and the fact that any normal open cover has a normal star refinement, it suffices to check that the following two diagrams commute.

\begin{diagram}
                                  &           &   \pi_1(\R(X,\mathcal U),x_0)       \\
                                  &    \ldTo  &                                     \\
\pi_1(\N(\St\mathcal U),\St U_0)  &           &    \dTo                             \\
                                  &    \rdTo  &                                     \\
                                  &           &   \pi_1(\R(X,\St\St\mathcal U),x_0) \\
\end{diagram}

\begin{diagram}
\pi_1(\N(\mathcal U),U_0)              &       &                                \\
                                       & \rdTo &                                \\
  \dTo                                 &       & \pi_1(\R(X,\St\mathcal U),x_0) \\
                                       & \ldTo &                                \\
\pi_1(\N(\St\St\mathcal U),\St\St U_0) &       &                                \\
\end{diagram}

Suppose $x_0,\ldots,x_n$ is a $\mathcal U$-chain in $X$ representing a loop in $\R(X,\mathcal U)$ based at $x_0$. Suppose $x_0,\ldots,x_n$ is sent to $\St U_0,\ldots,\St U_n$ which in turn is sent to $y_0,\ldots,y_n$. Now $y_0=y_n=x_0$ by assumption. To see that $y_0,\ldots,y_n$ is $\St\St\mathcal U$-homotopic to $x_0,\ldots,x_n$, notice that for any $i<n$,  $x_i,y_i,x_{i+1},y_{i+1}\in\St\St U$ where $U\in \mathcal U$ contains $x_i$ and $x_{i+1}$.

Now suppose $U_0,\ldots,U_n$ is a sequence of vertices in $\N(\mathcal U)$ that represents a loop in $\N(\mathcal U)$ based at $U_0$. Suppose $U_0,\ldots,U_n$ is sent to $x_0,\ldots,x_n$ which in turn is sent to $\St\St V_0,\ldots,\St\St V_n$ where $V_i\in\mathcal U$. Now $V_0=V_n=U_0$ by assumption. To see that the loop represented by $\St\St V_0,\ldots,\St\St V_n$ is homotopic to the one represented by $U_0,\ldots,U_n$ in $\N(\St\St\mathcal U)$, notice that for any $i<n$, $x\in \St\St U_i\cap\St\St V_i\cap\St\St U_{i+1}\cap\St\St V_{i+1}$ where $x$ is an element in $U_i\cap U_{i+1}$.
\endproof

There is a natural homomorphism $\pi_1(X,x_0)\to \check\pi_1(X,x_0)$ from the fundamental group to the shape group. Suppose $\alpha$ is a path in $X$ and $\mathcal U$ is an open cover of $X$. Choose $\delta>0$ so that any subinterval of $[0,1]$ of length $\delta$ is sent by $\alpha$ to some $U\in \mathcal U$. Define a $\mathcal U$-chain $\varphi_{\mathcal U}(\alpha)=\alpha(0),\alpha(\delta),\alpha(2\delta),\ldots,\alpha(1)$. This definition is independent of the choice of $\delta$; given $\delta_1<\delta_2$, the corresponding $\mathcal U$-chains will be $\mathcal U$-homotopic. Simply add the two chains together according to the order on $[0,1]$ to get another $\mathcal U$-chain which is $\mathcal U$-homotopic to both.

A similar argument shows that if $\mathcal V$ refines $\mathcal U$, $\varphi_{\mathcal V}(\alpha)$ will be $\mathcal U$-homotopic to $\varphi_{\mathcal U}(\alpha)$. Thus we have a generalized path $\varphi(\alpha)=\{\varphi_{\mathcal U}(\alpha)\}$. We show that $\varphi$ induces a well-defined homomorphism $\pi_1(X,x_0)\to \check\pi_1(X,x_0)$. Suppose $\alpha$ is fixed endpoint homotopic to $\beta$. Let $\mathcal U$ be a cover of $X$ and $\delta>0$ be such that any square $I\times I\subset [0,1]\times[0,1]$ of side length $\delta$ is sent by the homotopy to an element of $\mathcal U$. We have a sequence of paths $\alpha_{0},\alpha_{\delta},\alpha_{2\delta},\ldots,\beta$ given by the homotopy. To see that $\varphi_{\mathcal U}(\alpha_{i\delta})$ is $\mathcal U$-homotopic to $\varphi_{\mathcal U}(\alpha_{(i+1)\delta})$, notice the following chain is $\mathcal U$-homotopic to both.

\begin{center}
$\alpha_{i\delta}(0)\ \alpha_{i\delta}(\delta)\ \alpha_{(i+1)\delta}(\delta)\  \alpha_{(i+1)\delta}(2\delta)\ \alpha_{i\delta}(2\delta)\  \alpha_{i\delta}(3\delta)\ \alpha_{(i+1)\delta}(3\delta)\cdots$
\end{center}

We end the section by showing that this homomorphism is identical to the classical homomorphism. Given a cover $\mathcal U$, there is a map $\phi:X\to \N(\mathcal U)$ given by a partition of unity subordinate to $\mathcal U$, $\phi(x)=\sum\phi_U(x)U$. This map enjoys the property that $\phi^{-1}(\St U)\subset U$ where $\St U$ is the open star of the vertex $U$ in $\N(\mathcal U)$. $\phi$ induces a homomorphism $\pi_1(X)\to \check\pi_1(X)$. We show the following diagram commutes.

\begin{diagram}
             &       &   \pi_1(\R(X,\mathcal U),U_0)    \\
             & \ruTo &                                  \\
\pi_1(X,x_0) &       &              \dTo                \\
             & \rdTo &                                  \\
             &       & \pi_1(\N(\St\mathcal U),\St U_0) \\
\end{diagram}

Let $\phi:X\to \N(\St\mathcal U)$ be a map given by a partition of unity subordinate to $\St\mathcal U$. Let $\alpha$ be a path in $X$. Let $\delta>0$ so that if $I$ is a subinterval of $[0,1]$ of length $\delta$, $\alpha(I)$ is sent by $\phi$ to some $\St\St U$, the open star of a vertex $\St U$ of $\N(\St \mathcal U)$. Let $\delta$ also be such that $\alpha(I)$ lies in an element of $\mathcal U$.

Now $\alpha$ is sent to the chain $\alpha(0),\alpha(\delta),\alpha(2\delta),\ldots$ in $\R(X,\mathcal U)$ which in turn is sent to the chain $\St U_{0},\St U_{\delta},\St U_{2\delta},\ldots$ in $\N(\St\mathcal U)$ where $\alpha(k\delta)\in U_{k\delta}$ for $k\geq 0$. On the other hand, $\alpha$ is sent to the path $\phi\alpha$ in $\N(\St\mathcal U)$. We need to show that $\phi\alpha$ is fixed endpoint homotopic to the concatenation of edge paths associated with the chain $\St U_{0},\St U_{\delta},\St U_{2\delta},\ldots$. We proceed by induction on the number of terms in this chain. 

Now $\phi\alpha[0,\delta]\subset \St\St U$ and $\phi\alpha[\delta,2\delta]\subset\St\St U_1$ for some $U,U_1\in\mathcal U$.  Since $\phi\alpha(\delta)\in\St\St U\cap\St\St U_1$, it is in an open simplex having $\St U$ and $\St U_1$ as vertices. Now $\alpha[0,\delta]\subset\St U$ and $\alpha[\delta,2\delta]\subset\St U_1$ so $\alpha(\delta)\in\St U\cap\St U_1$ and $[\St U,\St U_1,\St U_{\delta}]$ is a simplex. The open simplex $[\St U,\St U_1,\St U_{\delta}]$ is contained in $\St\St U\cap\St\St U_1$ so we can join $\phi\alpha(\delta)$ to $\St U_{\delta}$ by a path $\gamma$ with $\gamma[0,1)\subset \St\St U\cap\St\St U_1$. Then, since $\phi\alpha[0,\delta]\subset\St\St U$, we can find a homotopy from $\phi\alpha[0,\delta]$ to the edge path $[\St U_0,\St U_{\delta}]$ that fixes $\alpha(0)=\St U_0$ and follows $\gamma$ from $\phi\alpha(\delta)$ to $\St U_{\delta}$.

For each $i>1$, $\phi\alpha[i\delta,(i+1)\delta]\subset \St\St U_i$ for some $U_i\in\mathcal U$. Suppose that $\phi\alpha[0,k\delta]$ is homotopic to the concatenation of edge paths associated with $\St U_{0},\St U_{\delta},\ldots\St U_{k\delta}$ where the homotopy fixes $\alpha(0)=\St U_0$ and follows a path $\gamma_k$ from $\phi\alpha(k\delta)$ to $\St U_{k\delta}$ with $\gamma_k[0,1)\subset\St\St U_{k-1}\cap\St\St U_k$. We follow the same procedure as above to find a path $\gamma_{k+1}$ from $\phi\alpha((k+1)\delta)$ to $\St U_{(k+1)\delta}$ that is contained in $\St\St U_k\cap\St\St U_{k+1}$. Since $\phi\alpha[k\delta,(k+1)\delta]\subset\St\St U_k$, we can find a homotopy from $\phi\alpha[k\delta,(k+1)\delta]$ to the edge path $[\St U_{k\delta},\St U_{(k+1)\delta}]$ that follows $\gamma_k$ from $\phi\alpha(k\delta)$ to $\St U_{k\delta}$ and $\gamma_{k+1}$ from $\phi\alpha((k+1)\delta)$ to $\St U_{(k+1)\delta}$.

\section{Homotopical Hausdorff}\label{HomotopicHaus}

We recall a standard topology on $\widetilde X$, the set of fixed endpoint homotopy classes of paths in $X$ starting at some basepoint $x_0$.

\begin{definition}\label{BasicTopology}
Given $[\alpha]\in\widetilde X$ with terminal point $x$ and a neighborhood $U$ of $x$ in $X$, $B([\alpha],U)$ is the set of all $[\beta]\in\widetilde X$ such that $\alpha^{-1}\beta$ is fixed endpoint homotopic to a path in $U$. We will call the topology generated by these sets the whisker topology on $\widetilde X$ following \cite{LocPath}.
\end{definition}

This topology is used in Spanier \cite{Spa} and Munkres \cite{Mun} for the classic construction of covering spaces. It is equivalent to the quotient topology inherited from $(X,x_0)^{(I,0)}$ under the compact open topology for locally path connected and semilocally simply connected spaces \cite[Lemma 2.1]{FisZas}.

Investigations into the structure of $\widetilde X$ leads one to realize that it can fail to be Hausdorff. The harmonic archipelago in \cite{BogSie} is a standard example of a Hausdorff space whose space of path homotopy classes is not Hausdorff. This situation motivates the following definition (see \cite{CanCon}).

\begin{definition}
A space $X$ is homotopically Hausdorff if for each $x\in X$ and each essential loop $\gamma$ based at $x$, there is a neighborhood $U$ of $x$ such that $\gamma$ is not fixed endpoint homotopic to a path in $U$.
\end{definition}

Notice that a space $X$ is homotopically Hausdorff if and only if for all $x\in X$, $\cap \pi_1(U,x)=1$ where $U$ runs over all neighborhoods of $x$. Also, for a Hausdorff space $X$, $X$ is homotopically Hausdorff if and only if $\widetilde X$ is Hausdorff for any basepoint.

It is shown in \cite{FisZas} that if a space is shape injective then it is homotopically Hausdorff. The space $A$ in \cite{Conetal} is an example of a Peano continuum that is homotopically Hausdorff but not shape injective (see Example \ref{exA}).

Investigation into the structure of $\pi_1(X)$ as a subspace of $\widetilde X$ in \cite{LocPath} lead to the definition of a new topology on $\widetilde X$. The new topology is based on the following definition. Given an open cover $\mathcal U$, let $\pi(\mathcal U,x)$ be the subgroup of $\pi_1(X,x)$ generated by elements of the form $[\alpha\gamma\alpha^{-1}]$ where $\alpha$ is a path starting at $x$ and $\gamma$ is a loop in some $U\in\mathcal U$. These groups are used in Spanier \cite{Spa} to detect when a fibration with unique path lifting is a covering map. 

\begin{definition}\label{SmallLoop}
Given $[\alpha]\in\widetilde X$ with terminal point $x$, a normal open cover $\mathcal U$ of $X$, and a neighborhood $V$ of $x$ in $X$, $B([\alpha],\mathcal U,V)$ is the set of all $[\beta]\in\widetilde X$ such that $\alpha^{-1}\beta$ is fixed endpoint homotopic to a loop in $\pi(\mathcal U,x)$ concatenated with a path in $V$. We will call the topology generated by these basic sets the lasso topology.
\end{definition}

There are slight differences between the above definition and the one that appears in \cite{LocPath}. As in the definition of the shape group, here we restrict our attention to normal covers. Also, in \cite{LocPath} it is required that $V\in\mathcal U$. This requirement does not effect the topology generated by the sets.

We now define an analogous version of homotopical Hausdorff for the lasso topology.

\begin{definition}
A space $X$ is lasso homotopically Hausdorff if for each $x\in X$ and each essential loop $\gamma$ based at $x$, there is a normal open cover $\mathcal U$ of $X$ such that $[\gamma]\notin\pi_1(\mathcal U,x)$.
\end{definition}

Notice that a space $X$ is lasso homotopically Hausdorff if for all $x\in X$, $\cap \pi_1(\mathcal U,x)=1$ where $\mathcal U$ runs over all normal open covers of $X$. Also, for a Hausdorff space $X$, $X$ is lasso homotopically Hausdorff if and only if $\widetilde X$ is Hausdorff for any basepoint under the lasso topology.

This concept was investigated in \cite{FisZas} where it is shown that if $\cap\pi_1(\mathcal U)=1$ for some collection of open covers of $X$, then the endpoint map $\widetilde X\to X$ has unique path lifting (where $\widetilde X$ is given the whisker topology). 

\begin{example}\label{exA}
We show that the space $A$ from \cite{Conetal} is not lasso homotopically Hausdorff. Let $A'$ be the topologist's sine curve $\{(x,\sin(1/x)):x\in(0,1]\}\cup \{0\}\times [-1,1]$ rotated about its limiting arc. It is the \textit{surface} and \textit{central limit arc} portions of $A$. \textit{Connecting arcs} are added to form the Peano continuum $A$ and the authors show that if a loop in $A'$ is nullhomotopic in $A$ then it is nullhomotopic in $A'$ (Lemma 4.3). Choose a basepoint $x$ on the surface portion of $A$ and a loop $\beta$ that goes once around the surface portion. Since $\beta$ is essential in $A'$ it must be essential in $A$ as well. Given any neighborhood $U$ of a point on the central limit arc, $\beta$ is fixed endpoint homotopic to a loop of the form $\alpha\gamma\alpha^{-1}$ where $\gamma$ is contained in $U$. Thus $[\beta]\in\pi_1(\mathcal U,x)$ for all open covers $\mathcal U$ of $A$.  Note that $A$ is not shape injective.
\end{example}

We now see that this version of homotopical Hausdorff is equivalent to shape injectivity for Peano spaces. 

\begin{lemma}
Suppose $X$ is locally path connected. Any normal open cover $\mathcal U$ of $X$ has a normal open refinement composed of path connected sets.
\end{lemma}

\proof
Let $\mathcal U$ be an open cover with associated partition of unity $\{\phi_U\}$. Given $U\in\mathcal U$, decompose it into its path components $\{V_\alpha\}$. Since $X$ is locally path connected these components are open. Given $x\in X$, define $\psi_{V_\alpha}(x)=\phi_U(x)$ if $x\in V_{\alpha}$ and $\psi_{V_\alpha}(x)=0$ otherwise. Then $\{\psi_{V_\alpha}\}_{U\in\mathcal U}$ is a partition of unity subordinate to $\{V_{\alpha}\}_{U\in\mathcal U}$.
\endproof

\begin{theorem}\label{MainTheorem}
Suppose $X$ is path connected and locally path connected. Then $X$ is lasso homotopically Hausdorff if and only if it is shape injective.
\end{theorem}

\proof
The reverse direction is essentially \cite[Proposition 4.8]{FisZas}. We provide a proof here. Suppose $X$ is shape injective and that $[\beta]\in\cap\pi_1(\mathcal U,x)$. Given $\mathcal U$, let $\lambda$ be a path fixed endpoint homotopic to $\beta$ so that $\lambda=\lambda_1\cdots\lambda_n$ where each $\lambda_i=\alpha_i\gamma_i\alpha_i^{-1}$, $\alpha_i$ is a path starting at $x$, and $\gamma_i$ is a loop in some $U\in\mathcal U$. Send $\lambda$ to $\pi_1(R(X,\mathcal U),x)$. Then the image of $\gamma_i$ is $\mathcal U$-homotopic to the constant chain at the terminal point of $\alpha_i$ so the image of $\lambda_i$ is $\mathcal U$-homotopic to the constant chain at $x$. Thus the image of $\lambda$ is trivial and the image of $\lambda$ in $\check\pi_1(X,x)$ is trivial. Given that $X$ is shape injective, we have that $\lambda$ is trivial.

Now suppose that $X$ is lasso homotopically Hausdorff. Suppose $\beta$ is a loop in $X$ based at $x$ whose image in $\check\pi_1(X,x)$ is trivial. Given a cover $\mathcal U$, we wish to show $[\beta]\in\pi_1(\mathcal U,x)$. Let $\mathcal V$ be a cover so that $\St\mathcal V$ refines $\mathcal U$ and let $\mathcal W$ be a refinement of $\mathcal V$ composed of path connected sets. The image of $\beta$ in $R(X,\mathcal W)$ is $\mathcal W$-homotopic to the trivial chain at $x$. We proceed by induction on the number of steps in this simplicial homotopy.

The image of $\beta$ in $R(X,\mathcal W)$ is represented by a $\mathcal W$-chain $x_0,\ldots,x_n$, i.e.,  $\beta=\beta_0\cdots\beta_{n-1}$ where each $\beta_i$ is a path in some element of $\mathcal W$ from $x_i$ to $x_{i+1}$.

Suppose a step of the simplicial homotopy starts at the $\mathcal W$-chain $y_0,\ldots,y_m$. Suppose there is a $[\lambda]\in\pi_1(\mathcal U,x)$ and a path $\alpha=\alpha_0\cdots\alpha_{m-1}$ associated with $y_0,\ldots,y_m$ (that means each $\alpha_i$ is a path in some element of $\mathcal W$ from $y_i$ to $y_{i+1}$) such that $\beta$ is fixed endpoint homotopic to $\lambda\alpha$. We show the same thing can be said about the next chain in the simplicial homotopy.

Suppose the next step of the simplicial homotopy is obtained by vertex addition, say $\ldots,y_i,y_{i+1},\ldots$ to $\ldots,y_i,y,y_{i+1},\ldots$. Now $y_i,y,y_{i+1}\in W$ for some $W\in\mathcal W$. Join $y_i$ to $y$ by a path $\lambda_1$ in $W$ and join $y$ to $y_{i+1}$ by a path $\lambda_2$ in $W$. Let $\lambda=\alpha_0\cdots\alpha_{i+1}\lambda_2^{-1}\lambda_1^{-1}\alpha_i^{-1}\cdots\alpha_0^{-1}$ and $\alpha=\alpha_0\cdots\alpha_{i-1}\lambda_1\lambda_2\alpha_{i+1}\cdots\alpha_{m-1}$. Then  $\beta$ is fixed endpoint homotopic to $\lambda\alpha$, $[\lambda]\in\pi_1(\mathcal U,x)$, and $\alpha$ is associated with the new $\mathcal W$-chain.

Now suppose the next step of the simplicial homotopy is obtained by vertex deletion, say $\ldots,y_i,y_{i+1},y_{i+2},\ldots$ to $\ldots,y_i,y_{i+2},\ldots$. Then $y_i,y_{i+2}\in W$ for some $W\in\mathcal W$. Join $y_i$ to $y_{i+2}$ by a path $\lambda$ in $W$. Let $\lambda=\alpha_0\cdots\alpha_{i+1}\lambda^{-1}\alpha_{i-1}^{-1}\cdots\alpha_0^{-1}$ and $\alpha=\alpha_0\cdots\alpha_{i-1}\lambda\alpha_{i+2}\cdots\alpha_{m-1}$. Then $\beta$ is fixed endpoint homotopic to $\lambda\alpha$, $[\lambda]\in\pi_1(\mathcal U,x)$, and $\alpha$ is associated with the new $\mathcal W$-chain.

At the end of the simplicial homotopy the $\mathcal W$-chain is the trivial chain at $x$. Thus we have $\beta$ is fixed endpoint homotopic to $\lambda\alpha$ where $[\lambda]\in\pi_1(\mathcal U,x)$ and $\alpha$ is associated with the trivial chain. Thus $\alpha$ is a loop based at $x$ in some element of $\mathcal W$ so $[\lambda\alpha]\in\pi_1(\mathcal U,x)$.
\endproof

\begin{remark}
The requirement of local path connectivity cannot be removed. In \cite[Remark 4.9]{FisZas} the authors give an example of a path connected space that is lasso homotopically Hausdorff but not shape injective. The space is related to the space $B$ in \cite{Conetal} (see Example \ref{exB}). In fact the authors in \cite{FisZas} cite \cite{RepZas} which became \cite{Conetal}. 
\end{remark}

\begin{example}\label{exB}
Let $B'$ be the topologist's sine curve $\{(x,\sin(1/(1-x))):x\in(0,1]\}\cup \{1\}\times [-1,1]$ rotated about a vertical axis at the point $(0,\sin(1))$. This space is lasso homotopically Hausdorff (it is locally simply connected). Connecting arcs are used in \cite{Conetal} to obtain the Peano continuum $B$. Since $B$ is not shape injective it cannot be lasso homotopically Hausdorff. A loop that goes around the outer annulus can be pulled in to the surface portion creating lassos. 
\end{example}

\end{document}